# Eikonal Slant Helices and Eikonal Darboux Helices In 3-Dimensional Riemannian Manifolds

**Mehmet Önder[a], Evren Zıplar[b], Onur Kaya[a]**
[a]*Celal Bayar University, Faculty of Arts and Sciences, Department of Mathematics, Muradiye Campus, 45047 Muradiye, Manisa, Turkey.*
E-mails: mehmet.onder@cbu.edu.tr, onur.kaya@cbu.edu.tr
[b]*Çankırı Karatekin University, Faculty of Science, Department of Mathematics, Çankırı, Turkey*
E-mail: evrenziplar@karatekin.edu.tr

**Abstract**

In this study, we give definitions and characterizations of eikonal slant helix curves, eikonal Darboux helices and non-normed eikonal Darboux helices in three dimensional Riemannian manifold $M^3$. We show that every eikonal slant helix is also an eikonal Darboux helix. Furthermore, we obtain that if the curve $\alpha$ is a non-normed eikonal Darboux helix, then $\alpha$ is an eikonal slant helix if and only if $\kappa^2 + \tau^2 = \text{constant}$, where $\kappa$ and $\tau$ are curvature and torsion of $\alpha$, respectively.

## 1. Introduction

In differential geometry of curves, the special curves whose curvatures satisfy some special conditions have an important role and more applications. One of the well-known of such curves is helix curve in the Euclidean 3-space $E^3$ which is defined by the property that the tangent line of the curve makes a constant angle with a fixed straight line (the axis of the general helix) [2]. Therefore, a general helix can be equivalently defined as one whose tangent indicatrix is a planar curve. Certainly, the helices in $E^n$ correspond with those whose unit tangent indicatrices are contained in hyperplanes. In 1802, Lancret stated a classical result for the helices and this result was first proved by B. de Saint Venant in 1845: *A necessary and sufficient condition that a curve to be a general helix is that the ratio of the first curvature to the second curvature be constant i.e., $\kappa/\tau$ is constant along the curve, where $\kappa$ and $\tau$ denote the first and second curvatures of the curve, respectively* [13].

Of course, there exist more special curves in the space. One of them is slant helix which first introduced by Izumiya and Takeuchi by the property that the normal lines of curve make a constant angle with a fixed direction in the Euclidean 3-space $E^3$ [7]. Later, Slant helices have been studied by some mathematicians and new kinds of these curves also have been introduced [1,5,8,9,10].

Moreover, some mathematicians have defined new special curve according to the Darboux vector of a space curve in $E^3$. First, the curve of constant precession has been introduced by Scofield [12]. Scofield has defined a curve of constant precession as follows:

**Definition 1.1.** ([12]) A unit speed curve of constant precession is defined by the property that its Darboux vector revolves about a fixed line in space with constant angle and constant speed. A curve of constant precession is characterized by having

$$\kappa(s) = w\sin(\mu(s)),$$
$$\tau(s) = w\cos(\mu(s)),$$

where $w > 0$ and $\mu$ are constants.

Later, Zıplar, Şenol and Yaylı have introduced the Darboux helix which is defined by the property that the Darboux vector of a space curve makes a constant angle with a fixed direction and they have given the characterizations of this new special curve [15].

Recently, new kinds of helices and slant helices have been introduced by Zıplar, Şenol and Yaylı by considering a space curve with a function $f : M \to \mathbb{R}$, where $M$ is a Riemannian manifold [14]. They have used eikonal function defined below and given $f$-eikonal helices and $f$-eikonal slant helices.

**Definition 1.2.** ([3]) Let $M$ be a Riemannian manifold and $g$ be its metric. For the function $f : M \to \mathbb{R}$, it is said that $f$ is eikonal if $\|\nabla f\|$ is constant, where $\nabla f$ is gradient of $f$, i.e., $df(X) = g(\nabla f, X)$.

Moreover, let us recall that $\|X\| = \sqrt{g(X, X)}$ for $X \in TM$, where $TM$ is the tangent bundle of $M$.

There exist many applications of $\nabla f$ in mathematical physics and geometry. For instance, if $f$ is non-constant on connected $M$, then the Riemannian condition $\|\nabla f\|^2 = 1$ is precisely the eikonal equation of geometrical optics. So, on a connected $M$, a non-constant real valued function $f$ is Riemannian if $f$ satisfies this eikonal equation. In the geometrical optical interpretation, the level sets of $f$ are interpreted as wave fronts. The characteristics of the eikonal equation (as a partial differential equation), are then the solutions of the gradient flow equation for $f$ (an ordinary differential equation), $x' = \nabla f$, which are geodesics of $M$ orthogonal to the level sets of $f$, and which are parameterized by arc length. These geodesics can be interpreted as light rays orthogonal to the wave fronts (See [4] for details).

***Theorem 1.1.*** *([6]) A complete Riemannian manifold $M$ without boundary admits a non-trivial affine function if and only if $M$ is isometric to Riemannian product $N \times \mathbb{R}$.*

***Lemma 1.1.*** *([11]) Let $(M, g)$ be a complete connected smooth Riemannian manifold. Then the followings are equivalent to each other.*

*i)* $f$ is an affine function.

*ii)* $f$ is smooth and its gradient vector field $\nabla f$ is parallel.

*iii)* $f$ is smooth and its Hessian $D^2 f$ vanishes everywhere.

*iv)* $f$ is smooth and $\nabla f$ is a Killing vector field with $\|\nabla f\| \equiv \text{constant}$.

In this paper, we introduce eikonal Darboux helix, non-normed eikonal Darboux helix and eikonal slant helix in 3-dimensional Riemannian manifold and give some characterizations of these new special curves.

## 2. Eikonal Slant Helices and Eikonal Darboux Helices

Let $\alpha = \alpha(s) : I \subset \mathbb{R} \to M^3$ be an immersed curve in 3-dimensional Riemannian manifold $M^3$. The unit tangent vector field of $\alpha$ will be denoted by $T$. Also, $\kappa > 0$ and $\tau$ will denote the curvature and torsion of $\alpha$, respectively. Therefore if $\{T, N, B\}$ is the Frenet frame of $\alpha$ and $\nabla$ is the Levi-Civita connection of $M^3$, then one can write the Frenet equations of $\alpha$ as

$$\nabla_T T = \kappa N$$
$$\nabla_T N = -\kappa T + \tau B$$
$$\nabla_T B = -\tau N$$

[2]. The vector $W = \tau T + \kappa B$ is called Darboux vector. Then for the Frenet formulae we have $\nabla_T T = W \times T$, $\nabla_T N = W \times N$, $\nabla_T B = W \times B$, where "×" shows the vector product in $M^3$.

**Definition 2.1.** Let $M^3$ be a 3-dimensional Riemannian manifold with the metric $g$ and let $\alpha(s)$ be a Frenet curve with the Frenet frame $\{T, N, B\}$ in $M^3$. Let $f : M^3 \to \mathbb{R}$ be a eikonal function along curve α, i.e. $\|\nabla f\| = \text{constant}$ along the curve $\alpha$. If the function $g(\nabla f, N)$ is non-zero constant along $\alpha$, then $\alpha$ is called an $f$-eikonal slant helix curve. And, $\nabla f$ is called the axis of the $f$-eikonal slant helix curve $\alpha$.

**Definition 2.2.** Let $M^3$ be a Riemannian manifold with the metric $g$. Let us assume that $\alpha$ be a curve in $M^3$ with Frenet frame $\{T, N, B\}$, non-zero curvatures $\kappa, \tau$ and Darboux vector $W = \tau T + \kappa B$. Also, let $f : M^3 \to \mathbb{R}$ be an eikonal function along $\alpha$. If the unit Darboux vector

$$W_0 = \frac{\tau}{\sqrt{\kappa^2 + \tau^2}} T + \frac{\kappa}{\sqrt{\kappa^2 + \tau^2}} B$$

of the curve $\alpha$ makes a constant angle $\varphi$ with the gradient of the function $f$, that is $g(W_0, \nabla f) = \cos\varphi$ is constant along $\alpha$, then the curve $\alpha$ is called $f$-eikonal Darboux helix.

Especially, if $g(W, \nabla f) = $ constant, then $\alpha$ is called a non-normed $f$-eikonal Darboux helix. Then, we have the following corollary.

***Corollary 2.1.*** *A non-normed $f$-eikonal Darboux helix is an $f$-eikonal Darboux helix if and only if $\kappa^2 + \tau^2$ is constant.*

**Exammple 2.1.** We consider the Riemannnian manifold $M^3 = \mathbb{R}^3$ with the Euclidean metric $g$. Let

$$f : M^3 \to \mathbb{R}$$
$$(x, y, z) \to f(x, y, z) = x + y^2 + z^2$$

be a function defined in $M^3$. Then, the curve

$$\alpha : I \subset \mathbb{R} \to M^3$$
$$s \to \alpha(s) = \left(\frac{s}{\sqrt{2}}, \cos\frac{s}{\sqrt{2}}, \sin\frac{s}{\sqrt{2}}\right)$$

is an $f$-eikonal slant helix curve in $M^3$.

Firstly, we will show that $f$ is an eikonal function along the curve $\alpha$. If we compute $\nabla f$, we find out $\nabla f$ as $\nabla f = (1, 2y, 2z)$. So, we get

$$\|\nabla f\| = \sqrt{1 + 4(y^2 + z^2)} \ .$$

And, if we compute $\|\nabla f\|$ along the curve $\alpha$, we find out

$$\|\nabla f\| = \sqrt{5} = \text{constant} \ .$$

That is, $f$ is an eikonal function along $\alpha$.

Now, we will show that the function $g(\nabla f, N)$ is non-zero constant along $\alpha$. Since

$$\nabla f = \left(1, 2\cos\frac{s}{\sqrt{2}}, 2\sin\frac{s}{\sqrt{2}}\right)$$

along $\alpha$ and

$$N = \left(0, -\cos\frac{s}{\sqrt{2}}, -\sin\frac{s}{\sqrt{2}}\right),$$

we obtain

$$g(\nabla f, N) = -2 = \text{constant}$$

along $\alpha$. Finally, $\alpha$ is an $f$-eikonal slant helix curve in $M^3$.

On the other hand, non-normed Darboux vector of $\alpha$ is

$$W = \left(\frac{\sqrt{2}}{2}, 0, 0\right),$$

and curvatures are $\kappa = \tau = \frac{1}{2}$. Then we obtain that

$$g(\nabla f, W) = \frac{\sqrt{2}}{2} = \text{constant}$$

along $\alpha$. So, $\alpha$ is an $f$-eikonal non-normed Darboux helix curve in $M^3$. Since $\kappa, \tau$ are constants, $\alpha$ is also an $f$-eikonal Darboux helix curve in $M^3$.

Now, we give some Theorems concerned with $f$-eikonal slant helix curve and $f$-eikonal Darboux helix. Whenever we write $M^3$, we will consider $M^3$ as a 3-dimensional Riemannian manifold with the metric $g$ and complete connected smooth without boundary and we will assume $M^3$ is isometric to a Riemannian product $N \times \mathbb{R}$.

**Theorem 2.1.** *Let us assume that $f : M^3 \to \mathbb{R}$ be a non-trivial affine function and $\alpha(s)$ be an f-eikonal slant helix curve in $M^3$. Then, the following properties are hold:*

*i) The function*

$$\frac{\kappa^2}{(\tau^2 + \kappa^2)^{3/2}} \left(\frac{\tau}{\kappa}\right)' \tag{1}$$

*is a real constant.*

*ii) The axis of $f$-eikonal slant helix is obtained as*

$$\nabla f = \frac{n\tau}{\sqrt{\tau^2 + \kappa^2}} T + \cos\theta N + \frac{n\kappa}{\sqrt{\tau^2 + \kappa^2}} B,$$

*where $n$ is a non-zero constant.*

***Proof.*** *i)* Since $\alpha$ is an $f$-eikonal slant helix, we have $g(\nabla f, N) = \cos\theta = \text{constant}$. So, there exist smooth functions $a_1 = a_1(s)$, $a_2 = a_2(s) = \cos\theta$ and $a_3 = a_3(s)$ of arc length $s$ such that

$$\nabla f = a_1 T + \cos\theta N + a_3 B, \tag{2}$$

where $\{T, N, B\}$ is a basis of $TM^3$ (tangent bundle of $M^3$).

Since $f$ is an affine function, by Lemma 1.1, $\nabla f$ is parallel in $M^3$, i.e., $\nabla_T \nabla f = 0$ along $\alpha$. Then, if we take the derivative in each part of (2) in the direction $T$ in $M^3$ and use the Frenet equations, we get

$$(T[a_1] - \kappa\cos\theta)T + (a_1\kappa - a_3\tau)N + (T[a_3] + \tau\cos\theta)B = 0, \tag{3}$$

Since $T[a_i] = a_i'(s)$, $(1 \leq i \leq 3)$ in (3) and the Frenet frame $\{T, N, B\}$ is linearly independent, we have

$$\begin{cases} a_1'(s) - \kappa\cos\theta = 0, \\ a_1\kappa - a_3\tau = 0, \\ a_3'(s) + \tau\cos\theta = 0. \end{cases} \tag{4}$$

From the second equation of the system (4) we obtain

$$a_1 = \left(\frac{\tau}{\kappa}\right) a_3. \tag{5}$$

Since $f$ is an eikonal function along $\alpha$, we have $\|\nabla f\|$ is constant. Then (2) and (5) give that

$$\left(\frac{\tau}{\kappa}\right)^2 a_3^2 + a_3^2 + \cos^2\theta = \text{constant}, \tag{6}$$

and from (6) we can write

$$\left(\left(\frac{\tau}{\kappa}\right)^2 + 1\right) a_3^2 = n^2, \tag{7}$$

where $n^2$ is a non-zero constant. From (7) we have

$$a_3 = \pm\frac{n}{\sqrt{\left(\frac{\tau}{\kappa}\right)^2 + 1}}. \tag{8}$$

By taking the derivative of (8) with respect to $s$ and using the third equation of the system (4), we get

$$\frac{\kappa^2}{\left(\tau^2+\kappa^2\right)^{3/2}}\left(\frac{\tau}{\kappa}\right)' = \text{constant} \tag{9}$$

which is desired function.

**ii)** By direct calculation from (5) and (8), we have

$$a_1 = \frac{n\tau}{\sqrt{\kappa^2+\tau^2}}, \quad a_3 = \frac{n\kappa}{\sqrt{\kappa^2+\tau^2}},$$

where $n$ is a non-zero constant. Then, the axis of f-eikonal slant helix is

$$\nabla f = \frac{n\tau}{\sqrt{\tau^2+\kappa^2}}T + \cos\theta N + \frac{n\kappa}{\sqrt{\tau^2+\kappa^2}}B \tag{10}$$

from (2).

The above Theorem has the following corollary.

***Corollary 2.2.*** *Let us assume that* $f : M^3 \to \mathbb{R}$ *be a non-trivial affine function and* $\alpha(s)$ *be an f-eikonal slant helix curve in* $M^3$ *with non-zero curvatures* $\kappa$ *and* $\tau$. *Then, the curvatures* $\kappa$ *and* $\tau$ *satisfy the following non-linear equation system:*

$$\left(\frac{n\tau}{\sqrt{\tau^2+\kappa^2}}\right)' - \mu\kappa = 0, \quad \left(\frac{n\kappa}{\sqrt{\tau^2+\kappa^2}}\right)' + \mu\tau = 0. \tag{11}$$

***Theorem 2.2.*** *Let* $f : M^3 \to \mathbb{R}$ *be a non-trivial affine function. Then, every f-eikonal slant helix in* $M^3$ *is also an f-eikonal Darboux helix in* $M^3$.

***Proof.*** Let $\alpha$ be f-eikonal slant helix in $M^3$. Then, from Theorem 2.1, the axis of $\alpha$ is

$$\nabla f = \frac{n\tau}{\sqrt{\tau^2+\kappa^2}}T + \cos\theta N + \frac{n\kappa}{\sqrt{\tau^2+\kappa^2}}B. \tag{12}$$

Considering the unit Darboux vector $W_0$, equality (12) can be written as follows

$$\nabla f = nW_0 + \cos\theta N, \tag{13}$$

which shows that $\nabla f$ lies on the plane spanned by $W_0$ and $N$. Since $n$ is a non-zero constant, from (13), we have $g(\nabla f, W_0) = n$ is constant along $\alpha$, i.e, $\alpha$ is an $f$-eikonal Darboux helix in $M^3$.

**Theorem 2.3.** *Let us assume that $f : M^3 \to \mathbb{R}$ be a non-trivial affine function and $\alpha$ be a non-normed $f$-eikonal Darboux helix. Then $\alpha$ is an $f$-eikonal slant helix if and only if $\kappa^2 + \tau^2 = \text{constant}$.*

**Proof.** Since $\alpha$ is a non-normed $f$-eikonal Darboux helix, we have $g(W, \nabla f) = \text{constant}$. On the other hand, there exist smooth functions $a_1 = a_1(s)$, $a_2 = a_2(s)$ and $a_3 = a_3(s)$ of arc length $s$ such that

$$\nabla f = a_1 T + a_2 N + a_3 B, \tag{14}$$

where $a_1$, $a_2$, $a_3$ are assumed non-zero and $\{T, N, B\}$ is a basis of $TM^3$. Since $f$ is an affine function, by Lemma 1.1., $\nabla f$ is parallel in $M^3$, i.e., $\nabla_T \nabla f = 0$ along $\alpha$. Then, if we take the derivative in each part of (14) in the direction $T$ in $M^3$ and use the Frenet equations, we get

$$(a_1' - a_2 \kappa)T + (a_1 \kappa + a_2' - a_3 \tau)N + (a_3' + a_2 \tau)B = 0, \tag{15}$$

where $a_i'(s) = T[a_i]$, $(1 \le i \le 3)$. Since the Frenet frame $\{T, N, B\}$ is linearly independent, we have

$$\begin{cases} a_1' - a_2 \kappa = 0, \\ a_1 \kappa + a_2' - a_3 \tau = 0, \\ a_3' + a_2 \tau = 0. \end{cases} \tag{16}$$

Equality $g(W, \nabla f) = \text{constant}$ gives that

$$a_1 \tau + a_3 \kappa = \text{constant}. \tag{17}$$

Differentiating (17) and using the first and third equations of system (16) we obtain

$$a_1 \tau' + a_3 \kappa' = 0 \tag{18}$$

From (18) and the second equation of system (16) it follows

$$a_2' = a_3 \frac{(\kappa^2 + \tau^2)'}{2\tau'} \tag{19}$$

Since $a_3 \ne 0$, from (19) we see that $a_2 = a_2(s)$ is constant if and only if $\kappa^2 + \tau^2 = \text{constant}$, i.e, $\alpha$ is an $f$-eikonal slant helix if and only if $\kappa^2 + \tau^2 = \text{constant}$.

From Theorem 2.3, Definition 1.1 and Corollary 2.1, we have the following corollaries.

**Corollary 2.3.** *Let $\alpha$ be a non-normed $f$-eikonal Darboux helix. Then $\alpha$ is an $f$-eikonal slant helix if and only if $\alpha$ is a curve of constant precession.*

***Corollary 2.4.*** *Let $\alpha$ be a non-normed $f$-eikonal Darboux helix. Then $\alpha$ is an $f$-eikonal slant helix if and only if $\alpha$ is an $f$-eikonal Darboux helix.*